\documentclass[11pt]{amsart}

\pdfoutput=1

\usepackage{amsmath,amsthm, amssymb, amsfonts}
\usepackage[mathscr]{eucal}
\usepackage{graphicx}
\usepackage[usenames,dvipsnames]{color}
\usepackage{subfigure}

\usepackage{tikz}
\usetikzlibrary{arrows,automata}

\usepackage[colorlinks=true, urlcolor=NavyBlue, linkcolor=NavyBlue, citecolor=NavyBlue, pdftitle={Torsion and Open Book Decompositions}, pdfauthor={John B. Etnyre, David Shea Vela--Vick}, pdfsubject={Transverse Invariants}, pdfkeywords={Legendrian Links, transverse Links, Heegaard Floer homology, 57M27, 57R58}]{hyperref}

\newtheorem{theorem}{Theorem}[section]
\newtheorem{lemma}[theorem]{Lemma}
\newtheorem{conj}[theorem]{Conjecture}

\theoremstyle{definition}

\newtheorem{remark}[theorem]{Remark}

\theoremstyle{remark}

\numberwithin{equation}{section}

\newcommand{\Z}{{\mathbb Z}}

\newcommand{\HFKH}{\widehat{\mathrm{HFK}}}
\newcommand{\HFKM}{\mathrm{HFK}^-}

\newcommand{\HFH}{\widehat{\mathrm{HF}}}

\newcommand{\SFH}{\mathrm{SFH}}

\newcommand{\SFHB}{\overline{\mathrm{SFH}}}
\newcommand{\cb}{\overline{c}}

\begin{document}

\thispagestyle{empty}

\title[Torsion and Open Book Decompositions]{Torsion and Open Book Decompositions}
\author{John B. Etnyre}
\address{School of Mathematics \\ Georgia Institute of Technology}
\email{etnyre@math.gatech.edu}
\urladdr{\href{http://www.math.gatech.edu/~etnyre}{http://www.math.upenn.edu/\~{}etnyre}}
\author{David Shea Vela--Vick}
\address{Department of Mathematics \\ Columbia University}
\email{shea@math.columbia.edu}
\urladdr{\href{http://www.math.columbia.edu/~shea}{http://www.math.columbia.edu/\~{}shea}}

\date{\today}
\maketitle

\begin{abstract}
We show that if $(B,\pi)$ is an open book decomposition of a contact 3--manifold $(Y,\xi)$, then the complement of the binding $B$ has no Giroux torsion. We also prove the sutured Heegaard-Floer $c$-bar invariant of the binding of an open book is non-zero. 
\end{abstract}


\section{Introduction} 
\label{sec:intro_back}

Interpreting and understanding contact geometric notions in terms of open book decompositions has been a central theme in the study of contact structures on 3--manifold ever since Giroux's fundamental breakthrough \cite{Gi2} equating contact structures on a 3--manifold (up to isotopy) with open book decompositions up to positive stabilizations (and isotopy). Two particularly noteworthy examples of this theme can be seen in identifying the Stein fillability of a contact structure with the existence of an associated open book with monodromy a composition of positive Dehn twists \cite{AkbulutOzbagci01, Gi2} and relating tight contact structures to right veering monodromies, \cite{HKM}. 

Another fundamental property of a contact structure is Giroux torsion. Recall that a contact manifold $(Y,\xi)$ is said to have {\it Giroux torsion $n$} if there exists a contact embedding $\phi : (T^2 \times I, \xi_{2n\pi}) \to (Y,\xi)$, where the contact structure $\xi_{2n\pi}$ on $T^2 \times I$ (thought of as $\mathbb{R}^2/\mathbb{Z} \times I$) is given by $\xi_{2n\pi}(s,t) = \ker(\cos(2n\pi t) dx + \sin(2n\pi t) dy)$.  

Currently Giroux torsion is the only known mechanism for a manifold to admit more than a finite number of tight contact structures. Thus it plays a central role in the course classification of contact structures on a 3--manifold. One would like to understand the relation between the Giroux torsion of a contact structure and properties of an open books supporting the contact structure. Currently we have the following reasonable conjecture concerning this relationship. 

\begin{conj}\label{conj:bn}
The binding number $bn(\xi)$ of a tight contact structure is bounded below by the Giroux torsion of $\xi.$
\end{conj}

Recall,  in \cite{EtnyreOzbagci08} the binding number of a contact structure was defined to be the minimal number of binding components for an open book supporting the contact structure among those open books supporting the contact structure having minimal genus pages. Since we currently have no way to ``see'' torsion from the perspective of supporting open book decompositions, and since we have no geometric understanding of what the binding number might be telling us about a contact structure, this is a particularly intriguing conjecture.

If Conjecture~\ref{conj:bn} is true, one expects there to be some interaction between the binding of an open book and the Giroux torsion of the associated contact structure. The simplest such interaction would be for the binding to somehow intersect the Giroux torsion, thus a very weak form of the above conjecture states that the complement of the binding of an open book for any contact structure has no Giroux torsion. 
 
The first progress on this conjecture occurred in \cite{Ve}, where the second author used invariants of Legendrian and transverse knots defined by Lisca, Ozsv\'ath, Stipsicz and Szab\'o in \cite{LOSS} to  show that if $(B,\pi)$ is an open book decomposition supporting $(Y,\xi)$ with connected binding $B$, then the complement of $B$ has no Giroux torsion.  We extend this to any open book decomposition, with no assumption on the number of binding components, by proving the following result. 

\begin{theorem}\label{thm:binding}
	Let $(B,\pi)$ be an open book decomposition of a contact 3--manifold $(Y,\xi)$, then the complement of the binding $B$ has no Giroux torsion.  In particular, $B$ must intersect each Giroux torsion layer in $(Y,\xi)$ nontrivially.
\end{theorem}

The proof involves a non-vanishing result for an invariant of transverse knots in a contact manifold. Following Stipsicz and V\'ertesi in \cite{StV} one can assign to a Legendrian or transverse knot (or link) $L \subset (Y,\xi)$ an invariant $\cb(L)$ which takes values in a certain sutured Floer homology group whose isomorphism type only depends on the topological type of $L$. We show this invariant never vanishes for the binding of an open book.

\begin{theorem}\label{thm:bindingnotvanish}
	Let $(B,\pi)$ be an open book decomposition of a contact 3--manifold $(Y,\xi)$, then $\cb(B)\not=0.$
\end{theorem}

Our proof of Theorem~\ref{thm:binding} now follows from a result of Ghiggini, Honda and Van Horn--Morris \cite{GHV} that in our language says any transverse link $L$ whose complement has Giroux torsion also has vanishing invariant, $\cb(L)=0.$

\subsection*{Acknowledgements} 
\label{sub:acknowledgements}
 The first author was partially supported by the NSF Grant DMS-0804820. The second author was partially supported by an NSF Postdoctoral Fellowship DMS-0902924.



\section{Background definitions and results} 
\label{sec:back}\label{sub:eh_bar}

We assume throughout familiarity with basic definitions and facts from 3--dimensional contact geometry, including open book decompositions, convex surface theory and Legendrian and transverse knot theory.  We also assume basic familiarity with Heegaard and sutured Heegaard-Floer homology and the contact invariants defined therein.  This material can be found in \cite{Et1}, \cite{Et2}, \cite{Ju1} and \cite{HKM2} respectively. 

Recall that to a (balanced) sutured manifold $(Y,\Gamma)$ one can associate the sutured Heegaard-Floer homology groups $\SFH(Y,\Gamma),$ \cite{Ju1}. In particular if $Y$ is a manifold with boundary and $\xi$ is a contact structure on $Y$ such that $\partial Y$ is convex, then the the dividing set $\Gamma_{\xi}$ on $\partial Y$ makes $Y$ into a balanced sutured manifold.  In  \cite{HKM2}, Honda, Kazez and Mati\'c defined an invariant of $\xi,$
\[
c(Y,\xi)\in \SFH (-Y, \Gamma_\xi),
\] 
that generalizes Ozsv\'ath and Szab\'o's Heegaard-Floer contact invariant on closed manifolds. 
We will make use of the following  gluing theorem of Honda, Kazez and Mati\'c for this invariant.

\begin{theorem}[Honda-Kazez-Mati\'c \cite{HKM3}]\label{thm:gluing}
	Let $(Y_1,\xi_1)$ and $(Y_2, \xi_2)$ be a compact contact 3--manifolds with convex boundary, and suppose that $(Y_1,\xi_1) \subset (Y_2,\xi_2)$ so that $m$ components of $Y_2-\mathrm{int } (Y_1)$  contain no boundary component of $Y_2.$ Then there exists a map
	\[
	\phi_{\xi_2-\xi_1} : \SFH(-Y_1,\Gamma_{\xi_1}) \to \SFH(-Y_2,\Gamma_{\xi_2}) \otimes V^{\otimes m},
	\] 
where $V$ is the two dimensional vector space over $\Z/2$: $\HFH(S^1\times S^2)$. Moreover, 
	\[
	\phi_{\xi_2-\xi_1}(c(Y_1,\xi_1))= c(Y_2,\xi_2)\otimes (x\otimes\cdots\otimes x),
	\] 
where $x$ is the contact invariant for the tight contact structure on $S^1\times S^2.$
\end{theorem}
\begin{remark}
The map in this theorem is only well defined up to sign when $\Z$-coefficients are used, but in this paper we will only consider $\Z/2$-coefficients and we can thus ignore the sign ambiguity.  
\end{remark}

We also make repeated use of a vanishing theorem of Ghiggini, Honda and Van Horn--Morris from \cite{GHV}.  Using Theorem~\ref{thm:gluing}, the proof of the main theorem of \cite{GHV} implies the following result. 

\begin{theorem}[Ghiggini-Honda-Van Horn--Morris \cite{GHV}]\label{thm:torsionvanish}
	If $(Y,\xi)$ is a contact manifold with positive Giroux torsion, then the contact invariant $c(Y,\xi)$ vanishes.
\end{theorem}


In \cite{HKM2}, Honda, Kazez and Mati\'c defined an invariant of Legendrian knots taking values in an appropriate sutured Floer homology group associated to a given Legendrian knot.  Simply put, if $L \subset (Y,\xi)$ is a Legendrian knot, then this invariant $c(L)$ is the sutured contact invariant of the complement of an open standard neighborhood of $L$.

A connection between the Legendrian invariants defined in \cite{HKM2} and \cite{LOSS} was explored by Stipsicz and V\'ertesi in \cite{StV}.  There, Stipsicz and V\'ertesi apply Theorem~\ref{thm:gluing} to map $c(L)$ to an intermediate invariant which we dub the {\it $c$-bar invariant} and denote $\cb(L)$.  Stipsicz and V\'ertesi show that the $c$-bar invariant maps to the Legendrian (hat) invariant from \cite{LOSS} under a natural identification of their ambient groups.  In this paper we ignore the connection with the Legendrian hat-invariant and focus our attention on the $c$-bar invariant. 

An elementary argument shows that $\cb(L)$ is a strictly weaker invariant than $c(L)$.  For instance, if $(Y,\xi)$ is a contact manifold with $c(Y,\xi) \neq 0$, then $c(L) \neq 0$ for any Legendrian $L \subset (Y,\xi)$.  On the other hand, $\cb(L)$ vanishes for any positively stabilized Legendrian knot $L$, regardless of the ambient contact manifold it lives in.  Forthcoming work of the authors \cite{EtV} explores the gap between $c(L)$ and $\cb(L)$ in detail.  There, we define a new ``sutured'' invariant that behaves much like the $\HFKM$ Legendrian invariant of \cite{LOSS}.

We now detail Stipsicz and V\'ertesi's construction of $\cb(L)$.  If $L \subset (Y,\xi)$ is a Legendrian knot, let $\nu(L) \subset (Y,\xi)$ denote an open standard neighborhood of $L$.  Focusing our attention on the complement of this standard neighborhood, we consider the contact manifold obtained by attaching a basic slice to the torus boundary of $(Y-\nu(L), \xi|_{Y-\nu(L)})$ so that the resulting dividing set on the new boundary consists of precisely two meridional (measured with respect to the original knot $L$) dividing curves.

\begin{figure}[htbp]
	\centering
	\subfigure[ ]{\label{fig:StV} 
	\includegraphics[scale=0.18]{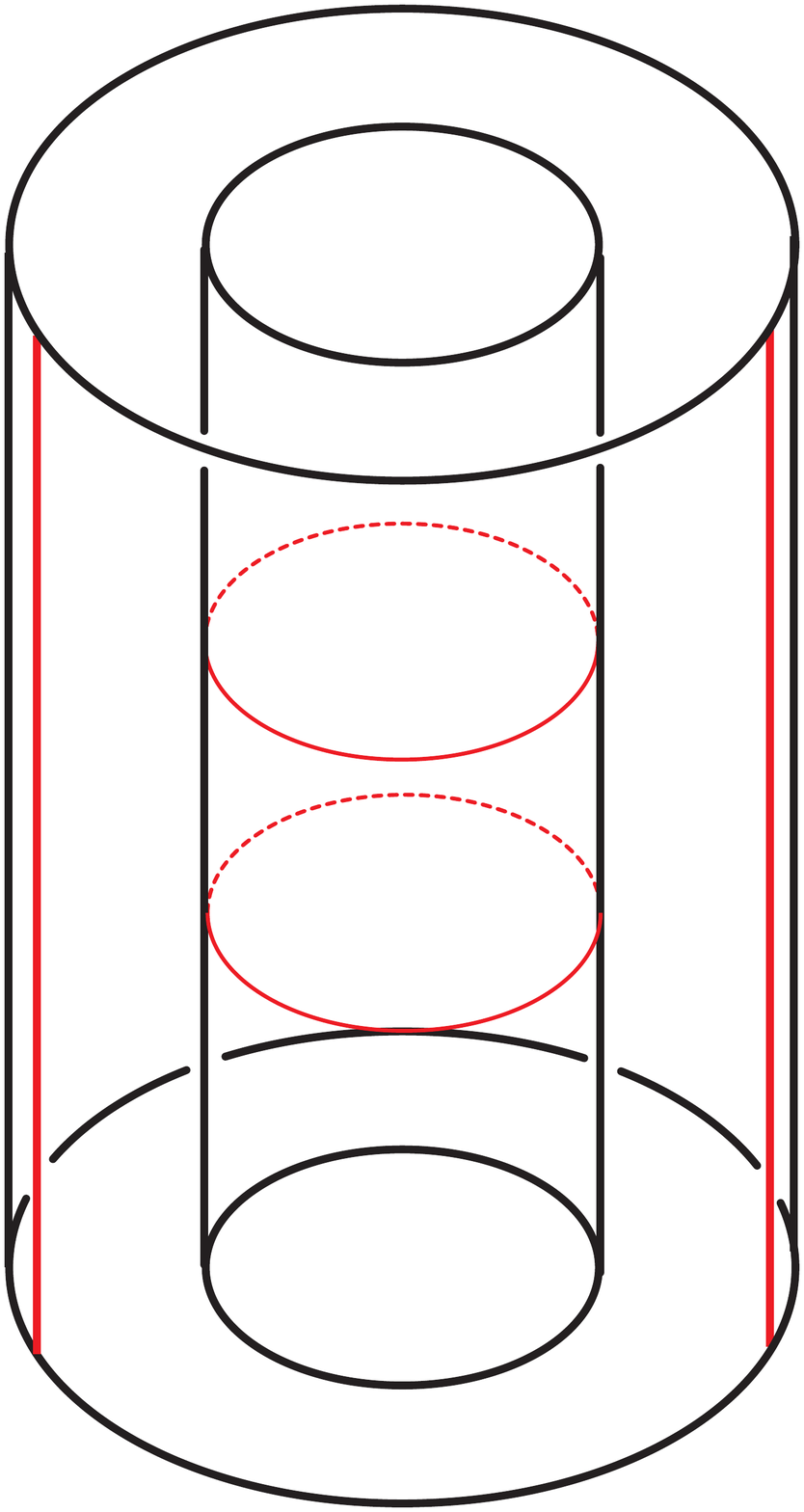}}
	\subfigure[ ]{\label{fig:StVfactor}
	\includegraphics[scale=0.18]{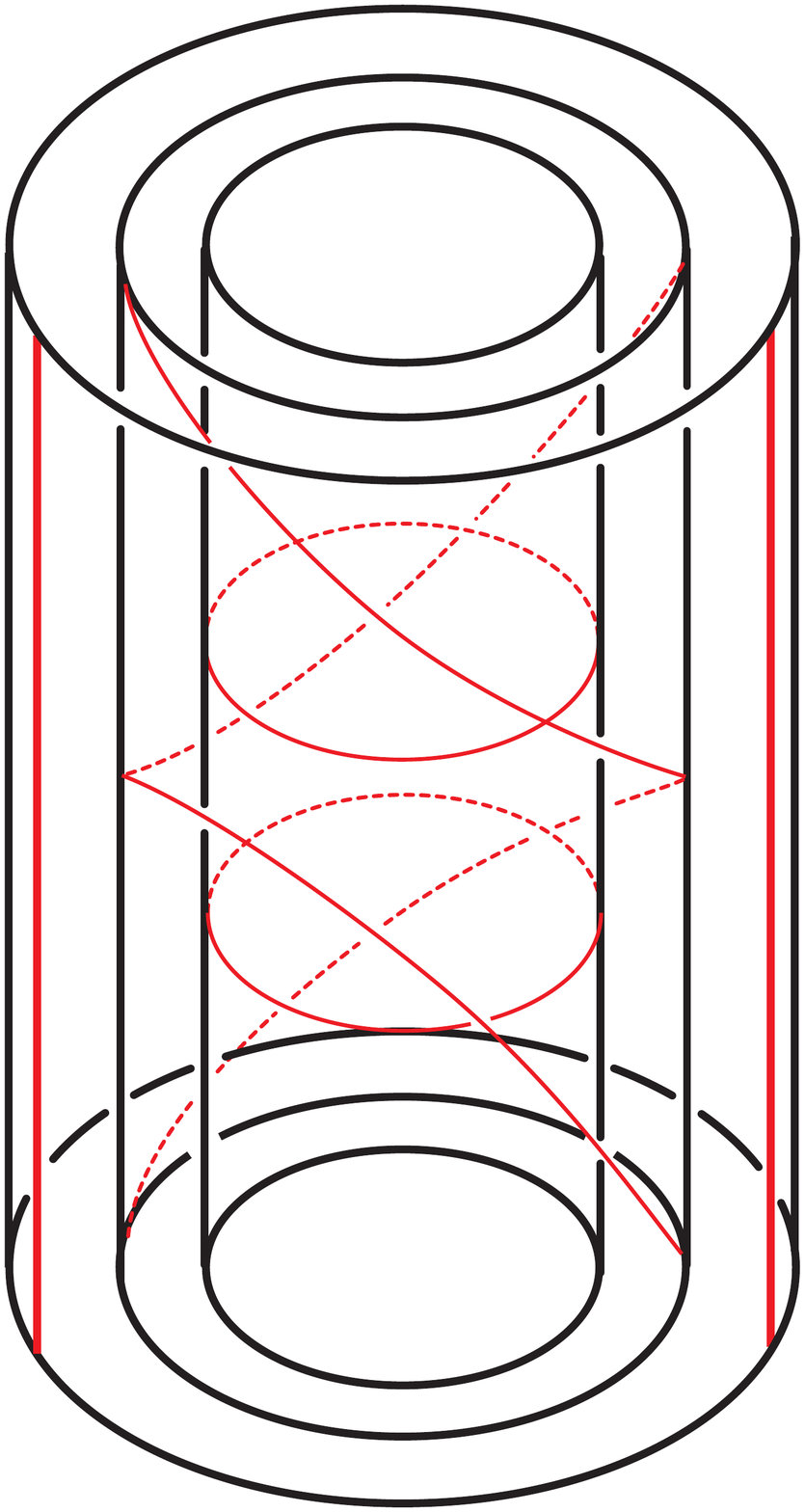}}
	\caption{Obtaining $\overline{Y}(L)$ and factor tori}
	\label{fig:HKMtoLoss}
\end{figure}

Denote this new contact manifold by $(\overline{Y}(L),\overline{\xi}_L)$.  A local picture of this construction is supplied in Figure~\ref{fig:StV}.  The contact invariant of the manifold $\overline{Y}(L)$ is the $c$-bar invariant of the Legendrian knot $L$.  We label this invariant by $\cb(L)$, and the sutured Floer homology of the ambient sutured manifold by $\SFHB(-Y,L)$.

Since the boundary sutures on $(\overline{Y}(L),\overline{\xi}_L)$ consist soley of meridians to $L$, the ambient sutured manifold, and therefore the group $\SFHB(-Y,L)$ depends only on the topological knot type of $L \subset Y$.  As a $\mathbb{Z}/2$-vector space, $\SFHB(-Y,L)$ is naturally isomorphic to the knot Floer homology group $\HFKH(-Y,L)$.

The invariant $\cb(L)$ is unchanged by negative stabilizations of the Legendrian knot $L$.  Recall that if $L$ and $L'$ differ by a negative stabilization, then the complements of open standard neighborhoods of $L$ and $L'$ differ by a basic slice attachment.  As depicted in Figure~\ref{fig:StVfactor}, the basic slice attachment yielding $\overline{Y}(L)$ factors into two basic slice attachments.  The first attachment corresponds to the stabilization $L \leadsto L'$, and the second corresponds to the attachment producing $\overline{Y}(L') = \overline{Y}(L)$.  Combined with an associativity property for the maps given by Theorem~\ref{thm:gluing}, see \cite{HKM3}, this proves invariance of $\cb(L)$ under negative stabilization.

Invariance of $\cb(L)$ under negative stabilizations, in turn, implies that $\cb$ induces an invariant of transverse knots: if $K \subset (Y,\xi)$ is a transverse knot, define $\cb(K) = \cb(L)$ for any Legendrian approximation $L$ of $K$.

This same construction extends to the case of Legendrian and transverse links, yielding invariants in this more general context.  Moreover, as observed by Stipsicz and V\'ertesi, the $c$-bar invariant is sensitive to Giroux torsion and satisfies the following vanishing theorem.

\begin{theorem}[Stipsicz-V\'ertesi \cite{StV}]\label{thm:torsion}
	If $L \subset (Y,\xi)$ be a Legendrian or transverse knot (or link) whose complement has positive Giroux torsion, then $\cb(L) = 0$.
\end{theorem}

\begin{proof}
If $L$ is a Legendrian link then the space $(\overline{Y}(L),\overline{\xi}_L)$ contains a copy of the complement $(Y-L, \xi|_{Y-L})$ as a proper subset.  Therefore, if the complement of $L$ has positive Giroux torsion, $(\overline{Y}(L),\overline{\xi}_L)$ must as well.  It follows immediately form the vanishing theorem in \cite{GHV}, Theorem~\ref{thm:torsionvanish}, and gluing theorem in \cite{HKM3}, Theorem~\ref{thm:gluing}, that $\cb(L) = 0$.

If $L$ is a transverse link with Giroux torsion in its complement, then there is a Legendrian approximation of $L$ that also has Giroux torsion in its complement.  Thus, it follows from the above argument that $\cb(L)=0.$
\end{proof}


\section{The $c$-bar invariant of bindings} 
\label{sec:ehb_bindings}

We prove Theorem~\ref{thm:binding} by showing that if $(B,\pi)$ is an open book decomposition of a given contact manifold $(Y,\xi)$, then $\cb(B)$ does not vanish.  It then follows from Theorem~\ref{thm:torsion} that the complement of $B$ is torsion free.

Throughout this section let $(B,\pi)$ be an open book decomposition for the contact manifold $(Y,\xi)$.  Choose a Legendrian approximation $L_B$ of $B$, and $S = \pi^{-1}(\theta_0)$ a fiber surface of $\pi$.

We now prove an auxiliary result that generalizes the universal tightness of the complement of the binding of an open book to the contact manifold $(\overline{Y}(L_B), \overline{\xi}_{L_B})$.

\begin{lemma}\label{lem:ybar_tight}
	The contact manifold $(\overline{Y}(L_B), \overline{\xi}_{L_B})$ is universally tight.
\end{lemma}

\begin{proof}
To show that $(\overline{Y}(L_B), \overline{\xi}_{L_B})$ is universally tight, we apply the Colin gluing theorem (see \cite{Co2}).
Colin's gluing theorem states that if $(Y_1,\xi_1)$ and $(Y_2, \xi_2)$ are two universally tight contact manifolds, and $T_i \subset Y_i$ are pre-Lagrangian, incompressible tori, then the contact manifold obtained by gluing $Y_1$ to $Y_2$ along $T_1$ and $T_2$ is universally tight.

Suppose for the moment that the binding $B$ has a single component, and let $T$ be the convex torus bounding $(Y - \nu(L_B), \xi|_{Y - \nu(L_B)})$.  The contact manifold $(\overline{Y}(L_B), \overline{\xi}_{L_B})$ is obtained from $(Y - \nu(L_B), \xi|_{Y - \nu(L_B)})$ by attaching a basic slices to the boundary torus $T$.

As observed in Section~\ref{sub:eh_bar}, this basic slice attachment can be factored as a composition of two basic slice attachments.  The first such attachment corresponds to a negative stabilization of $L_B$
\[
	(Y - \nu(L_B), \xi|_{Y - \nu(L_B)}) \; \leadsto \; (Y - \nu(L_B'), \xi|_{Y - \nu(L_B')}),
\]
while the second corresponds to the attachment,
\[
	(Y - \nu(L_B'), \xi|_{Y - \nu(L_B')}) \; \leadsto \; (\overline{Y}(L_B'), \overline{\xi}_{L_B'}) = (\overline{Y}(L_B), \overline{\xi}_{L_B}).
\]
Inside this first basic slice, we can find a pre-Lagrangian torus $T'$ parallel to the original boundary component $T$.  The complement of this pre-Lagrangian torus has two components; the first diffeomorphic to $T^2 \times I$, and the other to $\overline{Y - B}$.  

The contact structure restricted to either of these subspaces is universally tight.  More specifically, in the case of $T^2 \times I$, this is true because it sits as a ($\pi_1$-injective) subspace of a basic slice, a universally tight contact manifold.  Similarly, the second component is contained in complement of the binding, $(Y - B, \xi|_{Y-B})$, which is also universally tight. This fact is well-known but, as we could not find a direct proof of it in the literature, we give a proof here. 

By the definition of compatibility, there is a Reeb vector field $v$ for $\xi$ that is transverse to the pages of the open book and tangent to the binding. One can also arrange that $v$ is tangent to a family of concentric tori surrounding each of the binding components. Let $J$ be the standard almost complex structure associated to $v$ on the symplectization of $Y-B.$ That $v$ is tangent to concentric tori about $B$ implies that the corresponding ends of the symplectization are foliated by Levi-flat hypersurfaces, a fact which persists for finite covers of $Y-B.$ Thus, Hofer's proof that overtwisted contact structure on closed  3-manifold have contractible periodic Reeb orbits, \cite{Hofer93},  extends to the pair $(Y-B,v)$ and its finite covers.  As the Reeb vector field $v$ is transverse to the pages of the open book, $v$ can have no contractible periodic orbits.  Therefore, $\xi$ is tight, and remains tight when pulled back to finite covers. To see that the pullback of $\xi$ to the universal cover of $Y-B$ is tight, one can either use the fact that the fundamental group of $Y-B$ is residually finite, or that the universal cover of the symplectization of $Y-B$ has finite geometry at infinity.

In general, the binding $B$ consists of many components.  In this case, we apply the above argument to each component of $B$ yielding the contact manifold $(\overline{Y}(L_B), \overline{\xi}_{L_B})$.  All that remains to be checked before we can apply the Colin gluing theorem is that each of the boundary tori $T_i \subset \partial \, \overline{(Y - B)}$ are incompressible.

To see that each $T_i$ is incompressible, consider its preimage inside the universal cover of $Y - \nu(B)$.  We claim that each preimage is homeomorphic to a copy of $\mathbb{R}^2$, implying incompressibility.

Recall that $Y - \nu(B)$ is a surface bundle over $S^1$, fibered by oriented surfaces $S_\theta$, whose oriented boundary is $B$.  Consider first the intermediate cover that unwraps the $S^1$-factor.  This cover is homeomorphic to $S \times \mathbb{R}$.  In this intermediate cover, the preimage of each $T_i$ is a cylinder.

The universal cover of $Y - \nu(B)$ is homeomorphic to the universal cover of $S$ crossed with the real line $\mathbb{R}$.  Thus if $S$ is not a disk (our result in this case being obvious) the passing to the full universal cover, we see that the preimage of each of the above cylinders is homeomorphic to $\mathbb{R}^2$.

Hence, each $T_i$ is incompressible, completing the proof of Lemma \ref{lem:ybar_tight}. 
\end{proof}

\begin{lemma}\label{lem:nonvanish}
	The sutured contact invariant of $(\overline{Y}(L_B), \overline{\xi}_{L_B})$ does not vanish.
\end{lemma}

\begin{proof}
	As observed in \cite{Ve}, the Legendrian approximation $L_B$ can be chosen so that the twisting of the contact planes with respect to the framing induced on each component of $L_B$ by the fiber surface $S$ is $-1$.  In this case, the local picture around each boundary component of $(\overline{Y}(L_B), \overline{\xi}_{L_B})$ is shown in Figure~\ref{fig:YB}. In addition, as indicated in the figure, the dividing set on $S$ consists of one boundary parallel dividing curve for each boundary component of $S.$  This can be seen as follows.
	
	Let $B_1,\ldots, B_m$ be the components of $B.$ Stabilize the open book along $B_1$ (when stabilizing one chooses an arc in $S,$ take this arc to be a small boundary parallel arc). We get a new page $S'$ containing $S$ with $m+1$ boundary components $B'_1, B''_1, B_2,\ldots, B_m,$ where $B'_1$ is an unknot with self-lining $-1$ and $B''_1$ is transversely isotopic to $B_1.$  Moreover $\overline{S'-S}$ is a twice punctured disk with boundary $B_1,B_1', B_1'',$ see \cite{Ve}. Let $L_1,\ldots, L_m,$ be curves on $S'$ such that $L_i$ is isotopic to $B_i, i>1$ and $L_1$ isotopic to $B_1''.$ 
	
	As $S'$ is a page of an open book supporting $\xi$ we can make it convex (part of the definition of a contact structure being supported by an open book is that the Reeb vector field, which is a contact vector field, is transverse to the pages). Moreover, we can simultaneously Legendrian realize all of the $L_i$ on $S'.$ Note $B_i$ is the transverse push off of $L_i$ for all $i$ and the twisting of $L_i$ with respect to $S'$ is zero, for $i>1.$ 
	
	Let $S''$ be the subsurface of $S'$ with boundary $B_1, L_2,\ldots, L_m.$ We add an annulus to $S''$ along $B_1$ so that we get a surface, still denoted $S'',$ with boundary $L_1,L_2,\ldots, L_m.$ Moreover this surface is convex except on a disk touching the boundary component $L_1,$ (as it can be chosen to be a subsurface of $S'$ except along a disk touching the $L_1$ boundary component). The twisting of the contact planes along $L_1$ relative to $S''$ is $-1$ so we can isotop,  relative to where it was already convex, $S''$ so that it is convex. Thus the isotopy is supported near a disk $D$ touching $L_1.$ 
	
	The dividing set on $S''$ is empty except possibly in the disk $D$ since the Reeb field associated to $S'$ is transverse to $\xi.$ In the disk $D$ there must be a boundary parallel arc to account for the $-1$ twisting and nothing else (as the contact structure in the complement of the binding is tight). 
	
	Now negatively (Legendrian) stabilizing each $L_i, i>1,$ yields a new convex surface which we again denote $S$ with the desired dividing curves and the boundary Legendrian approximating $B.$

	\begin{figure}[htbp]
		\centering
		\subfigure[ ]{\label{fig:YB}
		\includegraphics[scale=0.2]{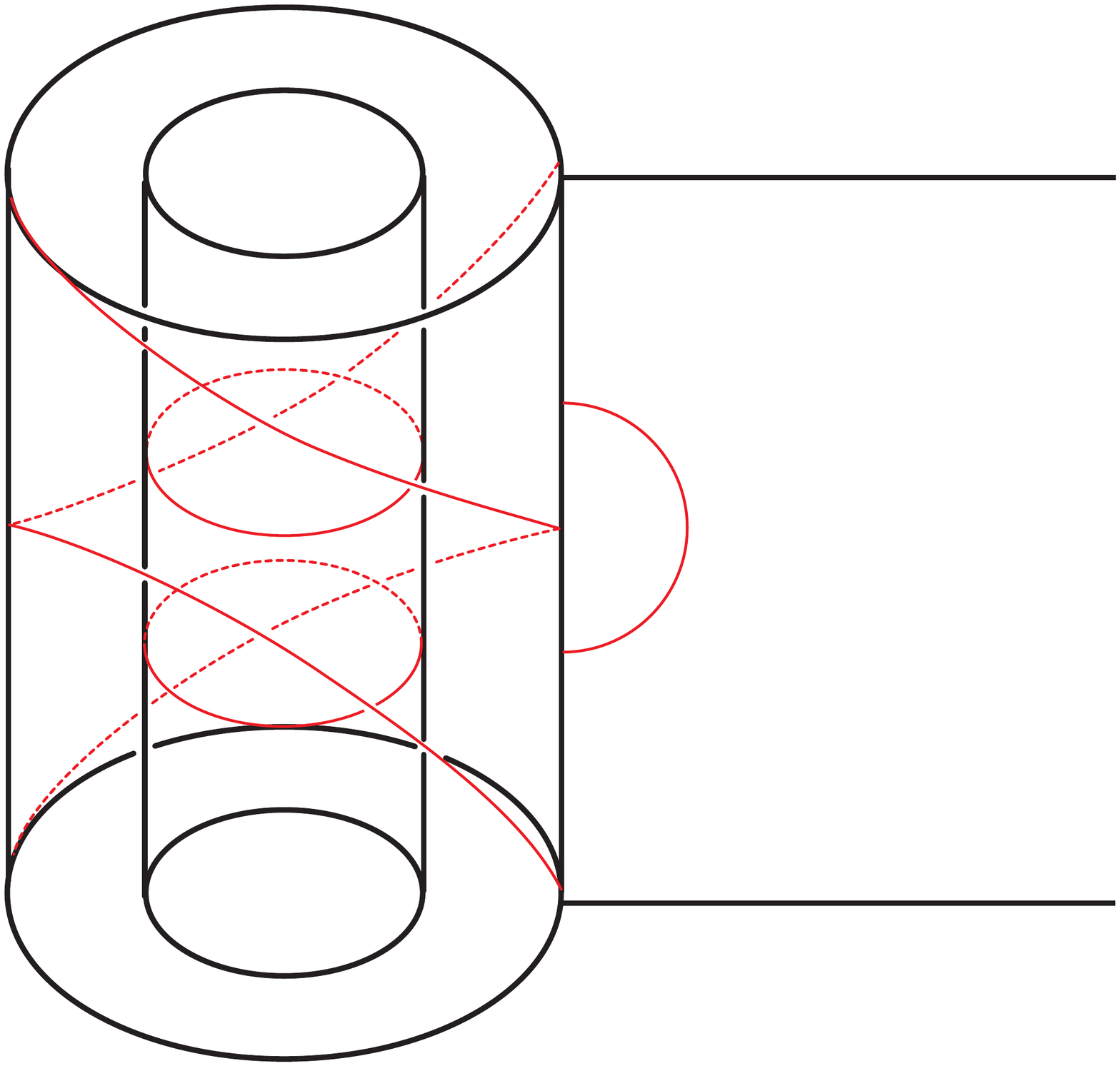}}
		\hspace{10px}
		\subfigure[ ]{\label{fig:basicslice}
		\includegraphics[scale=0.2]{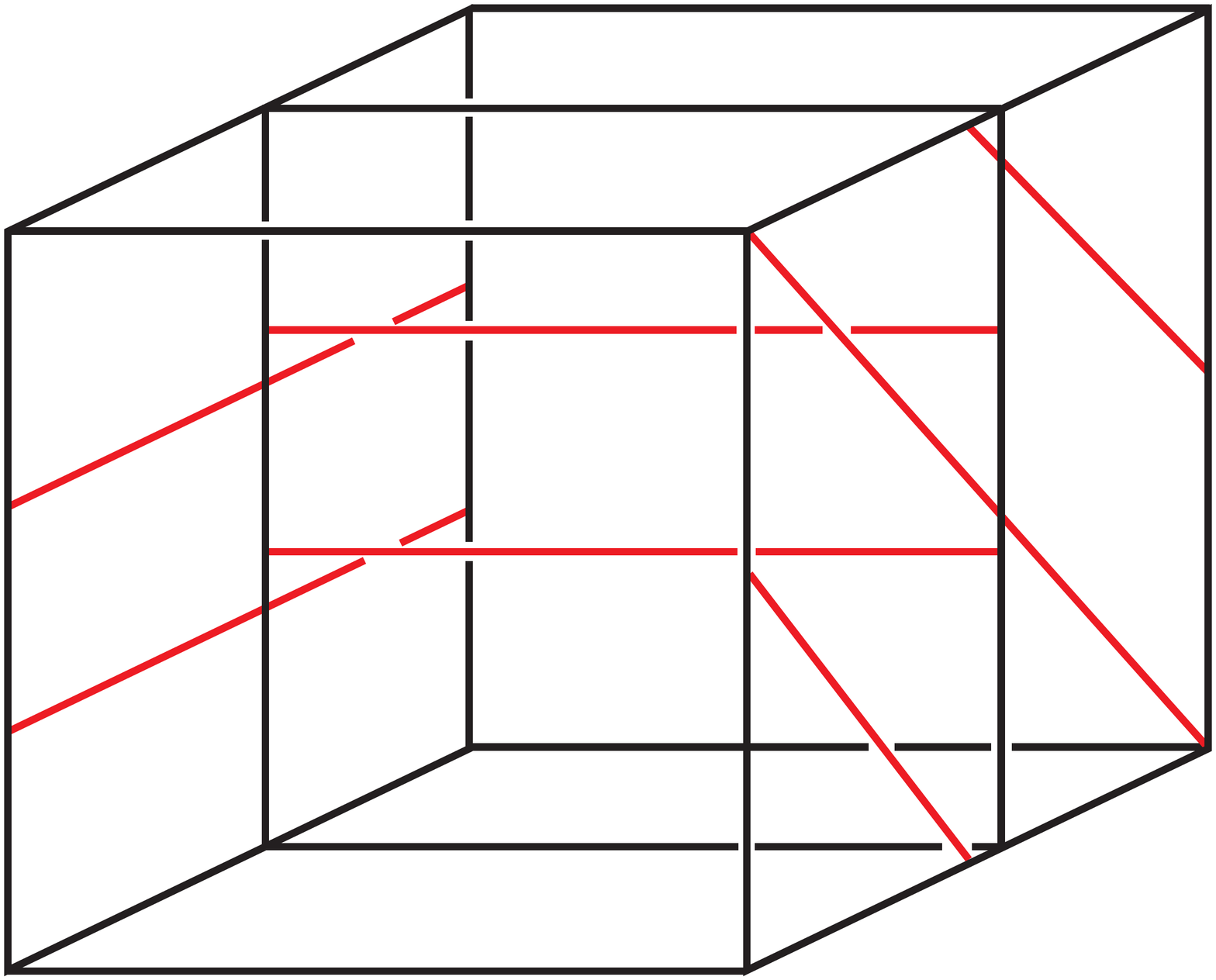}}
		\caption{Constructing $\overline{Y}(L_B)$ and viewing factor tori}
		\label{fig:page_and_slice}
	\end{figure}

	The dividing set on a convex annulus extending the page $S$ to the meridian-sloped boundary component consists, up to isotopy, of two horizontal dividing curves (see Figure~\ref{fig:basicslice}).  Denote by $S'$ the extension of the convex fiber surface $S$ by this convex annulus.

	The dividing set on $S'$ consists of a collection of boundary-parallel dividing curves; one for each boundary component of $\overline{Y}(L_B)$.  Such a surface is called {\it well-groomed}. It was shown in \cite{HKM2} that if $(Y_2,\xi_2)$ is obtained from $(Y_1,\xi_1)$ by cutting along a well-groomed convex surface, then the contact invariant of $(Y_2,\xi_2)$ is nonvanishing if and only if the contact invariant of $(Y_1,\xi_1)$ is nonvanishing.  Thus, it suffices to show that the contact invariant of the manifold obtained by cutting along $S'$ is nonzero.

	Cutting along $S'$, we obtain the tight contact manifold $(Y',\xi')$, whose dividing set near each of the original boundary tori is depicted in Figure~\ref{fig:YBcut}.
\begin{figure}[htbp]
		\centering
		\includegraphics[scale=0.3]{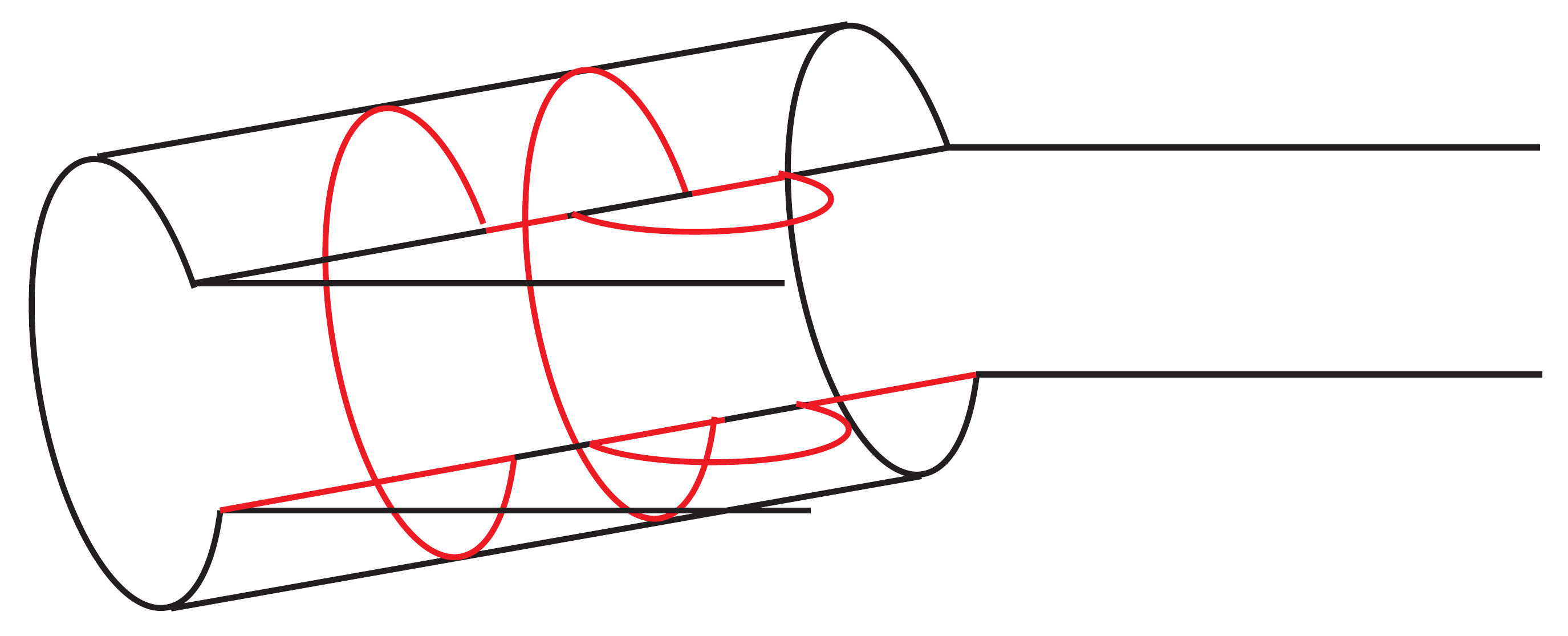}
		\caption{Splitting $\overline{Y}(L_B)$ along a page}
		\label{fig:YBcut}
	\end{figure}
	This contact manifold with convex boundary is tight, by Lemma~\ref{lem:ybar_tight}, and diffeomorphic to the surface $S' \times [0,1]$.  The dividing set is isotopic to the collection of curves $\partial S' \times {1/2}$.

Tight contact structures on such manifolds are unique, \cite{To}.  These contact structures are deformations of foliations, and embed into closed, Stein fillable contact manifolds. More specifically, notice that the dividing curves on the boundary of $S' \times [0,1]$ are the same as those on $S\times [0,1]$ (thought of as the complement of a standard neighborhood of $L_B$ in $(Y,\xi)$ cut open along the convex fiber surface $S$).  Thus the contact structures are contactomorphic. Moreover notice that if we take any take any open book with page $S$ and monodromy a composition of positive Dehn twists then the resulting contact manifold $(Y',\xi')$ is Stein fillable. Arguing as we just did for $S\times[0,1]$ and $S'\times[0,1]$, we see that if we remove a neighborhood of a Legendrian approximation of the binding from $Y'$ and cut the resulting manifold open along a page, the result will be a contact manifold contactomorphic to $S'\times[0,1].$ Thus $S'\times[0,1]$ contact embeds in a closed Stein fillable contact manifold. 

	As shown in \cite{OS4}, the contact invariant of a Stein fillable contact manifold is nonvanishing.  Therefore by the Honda-Kazez-Mati\'c gluing theorem \cite{HKM3}, Thoerem~\ref{thm:gluing}, it must also be the case that the sutured contact invariant of $(Y',\xi')$ is nonvanishing.
	\end{proof}

\begin{proof}[Proof of Theorem~\ref{thm:binding}]
Recall that, by definition, $\cb(B)$ is the contact invariant of $(\overline{Y}(L_B), \overline{\xi}_{L_B})$.  Therefore, by Lemma~\ref{lem:nonvanish}, $\cb(B) \neq 0$ for any $B$ which can be realized as the binding of an open book decomposition.  We now apply Theorem~\ref{thm:torsion} to conclude that the complement of such a $B$ has no Giroux torsion.
\end{proof}

\bibliographystyle{plain}
\nocite{*}
\bibliography{BindingNonvanish}

\begin{thebibliography}{10}

\bibitem{AkbulutOzbagci01}
Selman Akbulut and Burak Ozbagci.
\newblock Lefschetz fibrations on compact {S}tein surfaces.
\newblock {\em Geom. Topol.}, 5:319--334 (electronic), 2001.

\bibitem{Co2}
Vincent Colin.
\newblock Recollement de vari\'et\'es de contact tendues.
\newblock {\em Bull. Soc. Math. France}, 127(1):43--69, 1999.

\bibitem{Co1}
Vincent Colin.
\newblock Ouverts en g\'eom\'etrie de contact [d'apr\'es emmanuel giroux].
\newblock {\em S\'eminaire Bourbaki, 59 \`eme ann\'ee}, (969), 2006-2007.
\newblock to appear in Asterisque.

\bibitem{Et2}
John~B. Etnyre.
\newblock Legendrian and transversal knots.
\newblock In {\em Handbook of knot theory}, pages 105--185. Elsevier B. V.,
  Amsterdam, 2005.

\bibitem{Et1}
John~B. Etnyre.
\newblock Lectures on open book decompositions and contact structures.
\newblock In {\em Floer Homology, Gauge Theory, and Low-dimensional Topology},
  volume~5 of {\em Clay Math. Proc.}, pages 103--141. Amer. Math. Soc.,
  Providence, RI, 2006.

\bibitem{EH}
John~B. Etnyre and Ko~Honda.
\newblock Knots and contact geometry. {I}. {T}orus knots and the figure eight
  knot.
\newblock {\em J. Symplectic Geom.}, 1(1):63--120, 2001.

\bibitem{EtnyreOzbagci08}
John~B. Etnyre and Burak Ozbagci.
\newblock Invariants of contact structures from open books.
\newblock {\em Trans. Amer. Math. Soc.}, 360(6):3133--3151, 2008.

\bibitem{EtV}
John~B. Etnyre and David~Shea {Vela--Vick}.
\newblock Sutured {L}egendrian invariants and the ``{$U$}-action''.
\newblock {\em In preparation}, 2009.

\bibitem{GHV}
Paolo Ghiggini, Ko~Honda, and Jeremy {Van Horn--Morris}.
\newblock The vanishing of the contact invariant in the presence of torsion.
\newblock Preprint, \href{http://arxiv.org/abs/0706.1602v2}{\tt
  arXiv:0706.1602v2} [math.GT], 2007.

\bibitem{Gi1}
Emmanuel Giroux.
\newblock Structures de contact en dimension trois et bifurcations des
  feuilletages de surfaces.
\newblock {\em Invent. Math.}, 141(3):615--689, 2000.

\bibitem{Gi2}
Emmanuel Giroux.
\newblock G\'eom\'etrie de contact: de la dimension trois vers les dimensions
  sup\'erieures.
\newblock In {\em Proceedings of the International Congress of Mathematicians,
  Vol. II (Beijing, 2002)}, pages 405--414, Beijing, 2002. Higher Ed. Press.

\bibitem{Hofer93}
H.~Hofer.
\newblock Pseudoholomorphic curves in symplectizations with applications to the
  {W}einstein conjecture in dimension three.
\newblock {\em Invent. Math.}, 114(3):515--563, 1993.

\bibitem{Ho}
Ko~Honda.
\newblock On the classification of tight contact structures. {I}.
\newblock {\em Geom. Topol.}, 4:309--368 (electronic), 2000.

\bibitem{HKM1}
Ko~Honda, William Kazez, and Gordana Mati{\'c}.
\newblock On the contact class in {H}eegaard {F}loer homology.
\newblock Preprint, \href{http://arxiv.org/abs/math/0609734}{\tt
  arXiv:math/0609734v2} [math.GT], 2006.

\bibitem{HKM3}
Ko~Honda, William Kazez, and Gordana Mati{\'c}.
\newblock Contact structures, sutured {F}loer homology, and {TQFT}.
\newblock Preprint, \href{http://arxiv.org/abs/0807.2431}{\tt
  arXiv:math/0807.2431v1}[math.GT], 2008.

\bibitem{HKM}
Ko~Honda, William~H. Kazez, and Gordana Mati{\'c}.
\newblock Right-veering diffeomorphisms of compact surfaces with boundary.
\newblock {\em Invent. Math.}, 169(2):427--449, 2007.

\bibitem{HKM2}
Ko~Honda, William~H. Kazez, and Gordana Mati{\'c}.
\newblock The contact invariant in sutured {F}loer homology.
\newblock {\em Invent. Math.}, 176(3):637--676, 2009.

\bibitem{Ju1}
Andr{\'a}s Juh{\'a}sz.
\newblock Holomorphic discs and sutured manifolds.
\newblock {\em Algebr. Geom. Topol.}, 6:1429--1457 (electronic), 2006.

\bibitem{Ju2}
Andr{\'a}s Juh{\'a}sz.
\newblock Floer homology and surface decompositions.
\newblock {\em Geom. Topol.}, 12(1):299--350, 2008.

\bibitem{LOSS}
Paolo Lisca, Peter Ozsv{\'a}th, Andr{\'a}s~I. Stipsicz, and Zolt{\'a}n
  Szab{\'o}.
\newblock Heegaard {F}loer invariants of {L}egendrian knots in contact
  three-manifolds.
\newblock Preprint, \href{http://arxiv.org/abs/0802.0628}{\tt
  arXiv:0802.0628v1} [math.SG], 2008.

\bibitem{OS2}
Peter Ozsv{\'a}th and Zolt{\'a}n Szab{\'o}.
\newblock Holomorphic disks and three-manifold invariants: properties and
  applications.
\newblock {\em Ann. of Math. (2)}, 159(3):1159--1245, 2004.

\bibitem{OS1}
Peter Ozsv{\'a}th and Zolt{\'a}n Szab{\'o}.
\newblock Holomorphic disks and topological invariants for closed
  three-manifolds.
\newblock {\em Ann. of Math. (2)}, 159(3):1027--1158, 2004.

\bibitem{OS4}
Peter Ozsv{\'a}th and Zolt{\'a}n Szab{\'o}.
\newblock Heegaard {F}loer homology and contact structures.
\newblock {\em Duke Math. J.}, 129(1):39--61, 2005.

\bibitem{StV}
Andr{\'a}s~I. Stipsicz and Vera V{\'e}rtesi.
\newblock On invariants for {L}egendrian knots.
\newblock {\em Pacific J. Math.}, 239(1):157--177, 2009.

\bibitem{To}
Ichiro Torisu.
\newblock Convex contact structures and fibered links in 3-manifolds.
\newblock {\em Internat. Math. Res. Notices}, (9):441--454, 2000.

\bibitem{Ve}
David~Shea {Vela--Vick}.
\newblock On the transverse invariant for bindings of open books.
\newblock Preprint, \href{http://arxiv.org/abs/0806.1729}{\tt
  arXiv:0806.1729v1}[math.SG], 2008.

\end{thebibliography}

\end{document}